\documentclass[a4paper,11pt]{article}
\input epsf
\title{${\rm SL}_2(k)$ and a subset of words over $k$}
\author{Roland Bacher}
\begin{document}
\maketitle
%\par fichier su1.tex dans suitecurieuse. Pret pour soumission
%\bigskip
%version soumise le 11 septembre
%version finale avec corrections suggerees par les rapporteurs,
%renvoye le 4 mai
\par {\it Abstract:} 
\footnote{Math. class.: 05C38, 20G15, 68R15. Keywords: ${\rm SL}_2$,
${\rm PSL}_2$, equivalence relation, de Brujin sequence}
Given a field $k$, this paper defines a subset of the free semi-group
$F_k$ (whose elements are all finite words with letters in 
$k$) which has some interesting
arithmetic and combinatorial properties. The case $p=2$ has been treated in 
{\bf [B2]} with slightly different notations. 
The link with ${\rm SL}_2$ and motivation for studying this
subset originates in {\bf [B1]}.
\bigskip
\section{Introduction}
\medskip
Given a group or semi-group $\Gamma$ generated (in the sense of
semi-groups) by some subset $S\subset \Gamma$, the  
{\it Cayley graph} of $\Gamma$ with respect to $S$ is the 
oriented graph having vertices $\gamma\in \Gamma$ and oriented
edges $(\gamma,\gamma s)_{s\in S}$. 
In this paper, we consider generating sets $S$ which
are not necessarily finite thus yielding oriented Cayley
graphs which are perhaps not locally finite. 

The above situation gives rise to a homomorphism of semi-groups
$$\pi:F_S\longrightarrow \Gamma$$
where $F_S$ denotes the free semi-group on the set $S$
consisting of all finite words $s_1\dots s_l$ with letters 
in the alphabet $S$. The set $F_S$ is also called the {\it free
monoid} on $S$ and is often denoted by $S^*$. We will stick to the notation 
$F_S$ since our set $S$ will be identified with a field $k$ and 
the notation $k^*$ would be misleading in this case.

This paper deals with the group $\Gamma={\rm SL}_2(k)$ over an arbitrary field
$k$ generated by the set of matrices
$$S=\{\left(\begin{array}{cc} 0&-1\cr 1&\alpha\end{array}\right)\ \vert\ 
\alpha\in k\}.$$

We will study the preimage ${\cal A}\subset F_S$ of all finite words with
letters in $S$ (whose elements are indexed by the elements of the field $k$) such that
$$\pi({\cal A})=\{\left(\begin{array}{cc} a&-b\cr b^{-1}&0\end{array}\right)\ 
\vert \ a\in k\ ,\ b\in k^*\ \}\ .$$

The set ${\cal A}$ can also be described as follows: Recall that the projective
line ${\bf P}^1(k)$ consists of all $1-$dimensional subspaces in $k^2$.
We denote by $L(x)$ the subspace (line) spanned by $\left(\begin{array}{c}
1\cr x\end{array}\right)$ and by $L(\infty)$ the subspace 
spanned by $\left(\begin{array}{c} 0\cr 1\end{array}\right)$ (hence $L(x)$ denotes
the unique line in $k^2$ which has slope $x$ and runs through the origin).
The group ${\rm SL}_2(k)$ acts on ${\bf P}^1(k)$ (this action goes in fact
down to the projective group ${\rm PSL}_2(k)$ which is the quotient of
${\rm SL}_2(k)$ by $\pm \left(\begin{array}{cc}1&0\cr 0&1\end{array}\right)$).
The subset $\cal A$ considered above can now be defined as the set of all 
words $w\in F_S$ such that $\pi(w)L(\infty)=L(0)$. 

The aim of this paper is the description of some properties of the
set ${\cal A}\subset F_k$ and of its complement
${\cal C}=F_S\setminus {\cal A}\subset F_S$. All results will be stated in the next
section. Proofs will be given in section 3.

\section{Definitions and main results}

Consider a  field $k$ and the set 
$$S=\{\left(\begin{array}{cc} 0&-1\cr 1&\alpha\end{array}\right)\ \vert\
\alpha\in k\}\subset {\rm SL}_2(k)\quad .$$

{\bf Lemma 2.1.} {\sl The set $S$ generates ${\rm SL}_2(k)$ as a semigroup.}

This lemma shows that every element of ${\rm SL}_2(k)$ can be written
in at least one way as a finite word with letters in $S$. Since the
elements of $S$ are obviously indexed by $k$ we will only write
$\alpha$ instead of 
$\left(\begin{array}{cc} 0&-1\cr 1&\alpha\end{array}\right)$.
We identify hence the free monoid $F_S$ on $S$ with the free monoid
$F_k$ on $k$ consisting of the set of all
finite words with letters in the field $k$. We recall that we are interested in
the subset ${\cal A}\subset F_k$ defined by
$${\cal A}=\{w=\alpha_1\dots\alpha_l\in F_k\ \vert \ \pi(w)
L(\infty)=L(0)\}$$
(where $L(\infty)=k\left(\begin{array}{c}0\cr 1\end{array}\right)$ and 
$L(0)=k\left(\begin{array}{c}1\cr 0\end{array}\right)$)
and in its complement ${\cal C}=F_k\setminus {\cal A}$. 
We denote by $F_k^l$ the subset of all words of length
exactly $l$ in $F_k$. We set ${\cal A}^l={\cal A}\cap F_k^l$ and
${\cal C}^l={\cal C}\cap F_k^l$.

{\bf Theorem 2.2.} {\sl (i) $w\in {\cal A}^l$ if and only if 
$\alpha w\in {\cal C}^{l+1}$ and $w\alpha\in {\cal C}^{l+1}$ for every 
$\alpha\in k$.
\smallskip
\par \ \ (ii) For any $w\in {\cal C}^l$ there exist unique values
$\alpha,\ \beta\in k$ such that $\alpha w,\ w\beta\in {\cal A}^{l+1}$. 
\smallskip
\par \ \ (iii) if $\alpha_1\dots \alpha_l\in{\cal A}^l$ then $\alpha_2\alpha_3\dots
\alpha_l$ and $\alpha_1\alpha_2\dots \alpha_{l-1}\in {\cal C}^{l-1}$. 
\smallskip
\par \ \ (iv) $\alpha_1\alpha_2\dots\alpha_{l-1}\alpha_l\in {\cal A}^l$ 
if and only if
$\alpha_l\alpha_{l-1}\dots\alpha_2\alpha_1\in {\cal A}^l$.
\smallskip
\par \ \ (v) $\alpha_1\dots \alpha_l\in {\cal A}^l$ if and only if
$(-\alpha_1)\dots(-\alpha_l)\in {\cal A}^l$.}

{\bf Corollary 2.3.} {\sl For $k$ the finite field on $q=p^d$ elements
(with $p$ the finite prime characteristic of $k$) 
and for $l=0,1,2,\dots$ we have
$$\sharp({\cal A}^l)=\frac{q^l-(-1)^l}{q+1}\ ,\quad \sharp({\cal C}^l)=
\frac{q^{l+1}+(-1)^l}{q+1}\quad .$$
}

\par Consider the equivalence relation $\sim$ on $F_k$ with classes $\cal A$ and
$\cal C$. Denote by $\epsilon$ the empty word 
(of length $0$) in $F_k$. 
Extend the applications $x\longmapsto x+1,\ 
x\longmapsto x-1$ of the field $k$ to applications of the set
$k\cup \{\epsilon\}$ into itself by setting $(\epsilon\pm 1)=\epsilon$.

{\bf Proposition 2.4.} {\sl (i) One has 
$$\begin{array}{l}
x00y\sim xy\ ,\cr
x\alpha 1\beta y\sim x(\alpha-1)(\beta-1)y\ ,\cr
x\alpha (-1)\beta y\sim x(\alpha+1)(\beta+1)y\end{array}$$
where $x,y\in F_k,\ \alpha,\beta\in k\cup \{\epsilon\}$ with
$\alpha=\epsilon\Longrightarrow x=\epsilon$ and $\beta=\epsilon\Longrightarrow y=\epsilon$ (i.e. $\alpha$ is the last letter of word $x\alpha$ if $x\alpha$ is non-empty
and $\beta$ is the first letter of the word $\beta y$ if $\beta y$
is non-empty).

\ \ (ii) One has
$$\begin{array}{ll}
\alpha\beta x\sim (\beta-\alpha^{-1})x\quad &{\rm if}\  
0\not=\alpha\in k,\ \beta\in k,\ x\in F_k,\cr
0\beta x\sim x\quad& {\rm if}\ \beta\in k,\ x\in F_k\ .\end{array}$$
}

{\bf Remarks 2.5.} (i) If $k$ is the field on $2$ or $3$ elements,
then assertion (i) of previous proposition 
characterizes the sets $\cal A$ and $\cal C$ completely: 
it yields substitutions which
replace every word except $0\in {\cal A}$ and $\epsilon\in{\cal C}$ 
by an equivalent word which is strictly shorter.

\ \ (ii) Over any field $k$, assertion (ii) above and
the trivial observation $\alpha\sim\epsilon\Longleftrightarrow \alpha\not=0$
for $\alpha\in k$ determine the sets ${\cal A}$ and $\cal C$.

Set
$${\cal P}^l=\{\alpha_1\dots \alpha_l\in {\cal A}^l\ \vert\ 
\alpha_1\alpha_2\dots \alpha_h\in {\cal C}^h\quad {\rm for } \ h=1,\dots,
l-1\}$$
and ${\cal P}=\cup {\cal P}^l$.
\medskip
\par {\bf Theorem 2.6.} {\sl (i) (\lq\lq Unique factorization in ${\cal A}$'') 
We have $w\in {\cal A}$ if and only if $w$ can be written
as
$$w=p_1\delta_1p_2\delta_2\dots p_n\delta_n p_{n+1}$$
for some $n\geq 0$ with $p_1,\dots p_{n+1}\in {\cal P}$ and $\delta_1,\dots
,\delta_n\in k$. Moreover, such a factorization of $w\in {\cal A}$ is
unique.
\smallskip
\par \ \ (ii) We have for $l\geq 1$
$$\sharp({\cal P}^l)=(q-1)^{l-1}$$
if $k$ is the finite field on $q=p^d$ elements.
}

{\bf Corollary 2.7.} {\sl One has for any 
natural integer $l$ the identity
$$(x+1)\sum_{k=0}^{[l/2]} {l-k\choose k} x^k(x-1)^{l-2k}=x^{l+1}+(-1)^l$$
which is equivalent to the identities
$$\sum_{s=0}^k{l-s\choose s}{l-2s\choose k-s}(-1)^s=1$$
for $k=0,1,\dots,l$.}

{\bf Remark 2.8.} Theorem 2.6 shows that the vector space 
(over an arbitrary field) with basis the set
$$\{\epsilon\}\cup\{w\alpha\ \vert\ w\in {\cal A},\ \alpha\in k\}$$
can be turned into a graded algebra $\bf A$
(the product is given by extending linearly the concatenation of
words in $F_k$ and the grading is induced by the length of words in
$F_k$). It has in fact a very simple structure: the algebra $\bf A$ is
a free non-commutative algebra (on $q(q-1)^{l-2}$ generators of degree $l=2,3,4,\dots$ 
if $k$ is the finite field on $q=p^d$ elements).

Given two words $w,w'\in F_k^l$ of the form
$$w=\alpha_0\alpha_1\dots\alpha_{l-1},\ w'=\alpha_1\dots \alpha_{l-1}\alpha_l$$
we call $w'$ an {\it immediate successor} 
of $w$ and $w$ an {\it immediate predecessor} of $w'$.

{\bf Theorem 2.9.} {\sl Each element $w\in {\cal A}^l$ has a unique immediate
successor and a unique immediate predecessor in ${\cal A}^l$.}

Given an element $w_0 \in {\cal A}^l$, the previous theorem yields a sequence
$$w_0,\ w_1,\ w_2,\ w_3,\dots \in {\cal A}^l$$
with $w_{i+1}$ an immediate successor of $w_i$. 

Otherwise stated: For each $w\in {\cal A}^l$ there exists an infinite word
$$\tilde W=\dots \alpha_{-1}\alpha_0\alpha_1\alpha_2\alpha_3\dots$$
such that $\alpha_1\alpha_2\dots
\alpha_l=w$ and all factors $w_iw_{i+1}\dots w_{i+l-1}$
of length $l$ (subwords formed by $l$ consecutive letters)
of $\tilde W$ are elements of ${\cal A}^l$.

Until the end of this section we assume that $k$ is the finite field with
$q=p^d$ elements. In this case ${\cal A}^l$ is finite. Given $w\in
{\cal A}^l$
there exists hence a smallest integer $r$ such that the infinite word 
$\tilde W$ associated to $w$ is $r-$periodic.

{\bf Theorem 2.10.} {\sl Let $\tilde W=\dots \alpha_{r-1}\alpha_0\alpha_1\dots
\alpha_{r-1}\alpha_0\alpha_1\dots$ be an infinite 
$r-$periodic word with letters in a finite field $k$. 
Then there exists a smallest integer $t\leq q^2-1$ 
(in fact, $t$ is either $q$ or a divisor of $(q^2-1)$) such that all factors
of length $tr-1$ in $\tilde W$ belong to $\cal A$.}

{\bf Remark 2.11.} It follows (cf. assertion (i) in Lemma 3.1 of section 3)
that all factors of length $ltr-1$ ($l\geq 1$)
of $\tilde W$ belong also to $\cal A$. 
One can moreover show that if $m$ is an integer
with the property that all factors of length $m$ in $\tilde W$ belong to ${\cal A}^m$,
then $m=ltr-1$ for a suitable integer $l\geq 1$ (here $r$ denotes the minimal period
length of the infinite periodic word $\tilde W$).

{\bf Definition 2.12.} Given a finite set $E$ having $N\geq 2$ elements, a 
{\it mock parity check set} (MPCS for short) of length $d$ 
is a subset ${\cal M}\subset E^d$ (words of length $d$ with letters 
in $E$) such that
\par (i) each element $w\in {\cal M}$ has a unique immediate successor and a unique
immediate predecessor in ${\cal M}$.
\par (ii) ${\cal M}$ consists of exactly $N^{d-1}$ elements.

Denote by ${\rm Perm}_E$ the group of permutations  of the
finite set $E$ and let $\varphi:E^{d-2}
\longrightarrow {\rm Perm}_E$ be an application which 
associates to each element $z\in E^{d-2}$ a permutation 
$\varphi_z:E\longrightarrow E$.

{\bf Proposition 2.13.} {\sl The set
$${\cal M}=\{\alpha_1\alpha_2\dots\alpha_{d-1}\alpha_d\in E^d\ \vert\ 
\varphi_{\alpha_2\alpha_3\dots\alpha_{d-1}}(\alpha_1)=\alpha_d\}$$
is a MPCS and every MPCS is of this form.}

{\bf Remarks 2.14.} (i) This proposition shows that the set of all MPCS can be 
endowed with a group structure 
(the set of functions on $E^{N^{d-2}}$ with values in ${\rm Perm}_E$ has an
obvious group structure given by $(\varphi\psi)_z=\varphi_z\circ
\psi_z$). 
\smallskip
\par \ \ (ii) A MPCS ${\cal M}\subset
E^d$ yields a permutation of its elements: send each
$w\in{\cal M}$ to its (unique) successor in ${\cal M}$. 
Call a MPCS a (generalized) {\it de Brujin sequence} 
if the associated permutation consists of a unique cycle. One can
show that (generalized) de Brujin sequences exist for all integers $N\geq 2$ and 
$d\geq 1$.

{\bf Theorem 2.15.} {\sl Given a finite field $k$, the set
$${\cal M}^l={\cal A}^l\cup
\{\alpha_1\dots\alpha_l\in {\cal C}^l\ \vert\ \alpha_1\dots
\alpha_{l-1}\in {\cal A}^{l-1}\quad {\rm and}
\quad \alpha_2\dots\alpha_l\in {\cal A}^{l-1}\}$$
is a MPCS of $k^l$.}
%###############################################################################

\section{Proofs}

\par {\bf Proof of Lemma 2.1.} Take $\left(\begin{array}{cc} a&b\cr
c&d\end{array}\right) \in {\rm SL}_2(k)$. If $a\not=0$
we have necessarily $d=\frac{1+bc}{a}$ and the computation
$$\left(\begin{array}{cc}0&-1\cr 1&\frac{-c-1}{a} \end{array}\right)
\left(\begin{array}{cc}0&-1\cr 1&-a \end{array}\right)
\left(\begin{array}{cc}0&-1\cr 1&\frac{b-1}{a} \end{array}\right)
=\left(\begin{array}{cc} a&b\cr c&
\frac{1+bc}{a} \end{array}\right)$$
yields the result.

The case $a=0$ is reduced to precedent case by multiplying
first with $\left(\begin{array}{cc}0&-1\cr 1&0 \end{array}\right)$ and by remarking
that this matrix has order $4$ (or $2$ if the ground field $k$ is
of characteristic $2$).
\hfill QED

{\bf Proof of Theorem 2.2.} One has
$$\left(\begin{array}{cc} a&b\cr c&d\end{array}\right)
\left(\begin{array}{cc} 0&-1\cr 1&x\end{array}\right)=
\left(\begin{array}{cc} b&-a+bx\cr d&-c+dx\end{array}\right)$$
which shows that $w\alpha\not\in {\cal A}$ if $w\in \cal A$ (since then
$\pi(w)=\left(\begin{array}{cc}a&-b\cr b^{-1}&0\end{array}\right)$). 
On the other hand,
if $w\in {\cal C}$ then $\pi(w)=
\left(\begin{array}{cc} a&b\cr c&d\end{array}\right)$ 
with $d\not= 0$ and the above computation implies the existence of a unique 
$\beta$ such that $w\beta\in {\cal A}$. This proves half of (i) and (ii). The 
proof of the remaining half is similar (it is also implied by assertion (iv)).
\par In order to prove (iii) one considers
$$\begin{array}{c}
\displaystyle \pi(\alpha_2\dots \alpha_l)=
\pi(\alpha_1)^{-1}\pi(\alpha_1\dots\alpha_l)=
\left(\begin{array}{cc} \alpha_1&1\cr -1&0\end{array}\right)
\left(\begin{array}{cc} a&-b\cr b^{-1}&0\end{array}\right)\cr
\displaystyle =\left(\begin{array}{cc} \alpha_1 a+b^{-1}&-\alpha_1 b\cr
-a&b \end{array}\right)\end{array}$$
which shows that $\alpha_2\dots\alpha_l\in {\cal C}^{l-1}$. A similar computation
yields $\alpha_1\dots \alpha_{l-1}\in {\cal C}^{l-1}$.
\par Since 
$$\left(\begin{array}{cc} 0&1\cr 1&0\end{array}\right)
\left(\begin{array}{cc} 0&-1\cr 1&\alpha\end{array}\right)
\left(\begin{array}{cc} 0&1\cr 1&0\end{array}\right)
=\left(\begin{array}{cc} \alpha&1\cr -1&0\end{array}\right)$$
we get by conjugating $\pi(\alpha_1\dots\alpha_l)=
\left(\begin{array}{cc} a&b\cr c&d\end{array}\right)$ with $\sigma=
\left(\begin{array}{cc} 0&1\cr 1&0\end{array}\right)$
$$\sigma \pi(\alpha_1\dots\alpha_l)\sigma=
\left(\begin{array}{cc} \alpha_1&1\cr -1&0\end{array}\right)\cdots
\left(\begin{array}{cc} \alpha_l&1\cr -1&0\end{array}\right)
=\left(\begin{array}{cc} d&c\cr b&a\end{array}\right)$$
If $\pi(\alpha_1\dots\alpha_l)=
\left(\begin{array}{cc} a&-b\cr b^{-1}&0\end{array}\right)$ we get by taking the inverse of
$\sigma \pi(\alpha_1\dots\alpha_l)\sigma$
$$\pi(\alpha_l\dots\alpha_1)=
\left(\begin{array}{cc} 0&-1\cr 1&\alpha_l\end{array}\right)\cdots
\left(\begin{array}{cc} 0&-1\cr 1&\alpha_1\end{array}\right)
=\left(\begin{array}{cc} 0&b^{-1}\cr -b&a\end{array}\right)^{-1}
=\left(\begin{array}{cc} a&-b^{-1}\cr b&0\end{array}\right)$$
which shows that $\alpha_l\dots\alpha_1\in {\cal A}$ and proves (iv).
\par Transposing $\pi(\alpha_1\dots\alpha_l)$ and multiplying by $(-1)^l$ shows
that $(-\alpha_l)\dots(-\alpha_1)\in{\cal A}$. Assertion (iv) implies now (v).
\medskip
\par {\bf Remark 3.1.} The properties of the action
of ${\rm SL}_2(k)$ (or ${\rm PSL}_2(k)$) on the projective line ${\bf P}^1(k)$
can be used to get a more conceptual proof of most assertions in Theorem 2.2.
\medskip
\par {\bf Proof of Corollary 2.3.} Assertion (ii) of Theorem 2.2 shows that
$\sharp({\cal A}^{l+1})\geq \sharp ({\cal C}^l)$ and assertion (iii) implies
$\sharp({\cal A}^{l+1})\leq \sharp ({\cal C}^l)$ hence establishing
$\sharp({\cal A}^{l+1})= \sharp ({\cal C}^l)$. Induction on $l$ (using
the obvious identity $\sharp({\cal A}^l)+\sharp ({\cal C}^l)=q^l$) yields now the result.
\medskip
\par {\bf Proof of Proposition 2.4.} The first line of assertion (i)
follows from the identity
$$\left(\begin{array}{cc} 0&-1\cr 1&0 \end{array}\right)
\left(\begin{array}{cc} 0&-1\cr 1&0 \end{array}\right)=
\left(\begin{array}{cc} -1&0\cr 0&-1\end{array}\right)\ .$$
If $\alpha$ and $\beta$ are both
non-empty, the last two lines of the proposition follow from the identities
$$\left(\begin{array}{cc} 0&-1\cr 1&\alpha \end{array}\right)
\left(\begin{array}{cc} 0&-1\cr 1&1\end{array}\right)
\left(\begin{array}{cc} 0&-1\cr 1&\beta \end{array}\right)$$
$$=\left(\begin{array}{cc} -1&1-\beta\cr \alpha-1&\alpha\beta-\alpha-\beta 
\end{array}\right)=
\left(\begin{array}{cc} 0&-1\cr 1&\alpha-1 \end{array}\right)
\left(\begin{array}{cc} 0&-1\cr 1&\beta-1 \end{array}\right)$$
and
$$\left(\begin{array}{cc} 0&-1\cr 1&\alpha \end{array}\right)
\left(\begin{array}{cc} 0&-1\cr 1&-1\end{array}\right)
\left(\begin{array}{cc} 0&-1\cr 1&\beta \end{array}\right)=
\left(\begin{array}{cc} 1&1+\beta\cr -1-\alpha&-\alpha-\beta-\alpha\beta \end{array}\right)
$$
$$=-\left(\begin{array}{cc} 0&-1\cr 1&\alpha+1 \end{array}\right)
\left(\begin{array}{cc} 0&-1\cr 1&\beta+1 \end{array}\right)$$
(in fact, the last line of Proposition 2.4 is easily deduced from the second one by using 
assertion (v) of Theorem 2.2).
We leave the remaining cases of assertion (i)
(with $\epsilon\in\{\alpha,\beta\}$) to the reader (they follow also easily from Theorem 2.6).

Assertion (ii) follows from the computations
$$\left(\begin{array}{cc}0&-1\cr 1&\alpha\end{array}\right)
\left(\begin{array}{cc}0&-1\cr 1&\beta\end{array}\right)
\left(\begin{array}{cc}\alpha\beta a+\beta c-a&\alpha\beta b+\beta d-b\cr 
-\alpha a-c&-\alpha b-d\end{array}\right)=
\left(\begin{array}{cc}a&b\cr c&d\end{array}\right)$$
and
$$\left(\begin{array}{cc}0&-1\cr 1&\gamma\end{array}\right)
\left(\begin{array}{cc}\alpha\beta a+\beta c-a&\alpha\beta b+\beta d-b\cr 
-\alpha a-c&-\alpha b-d\end{array}\right)$$
$$=
\left(\begin{array}{cc}
\alpha a+c&\alpha b+d\cr \alpha\beta a+\beta c-a-\alpha\gamma a-\gamma c&
\alpha\beta b+\beta d-b-\alpha\gamma b-\gamma d\end{array}\right)\ .$$
\medskip
\par {\bf Lemma 3.2.} (i) If $w,w'\in {\cal A}$ then $ww'\in {\cal C}$ and 
$w\alpha w'\in {\cal A}$ for any $\alpha\in k$.
\smallskip
\par \ \ (ii) If exactly one of $w,w'$ is an element of ${\cal A}$ then
$w\alpha w'\in {\cal C}$ for any $\alpha\in k$.
\medskip
\par {\bf Proof of Lemma 3.2.} The computation
$$\left(\begin{array}{cc}a&-b\cr b^{-1}&0\end{array}\right)
\left(\begin{array}{cc}0&-1\cr 1&\alpha \end{array}\right)
\left(\begin{array}{cc}a'&-b'\cr b'^{-1}&0\end{array}\right)
=\left(\begin{array}{cc}-ba'-ab'^{-1}-\alpha bb'^{-1}&bb'\cr
-(bb')^{-1}&0 \end{array}\right)$$
shows (i).
\par Let us now suppose that $w\in {\cal A},\ w'\in {\cal C}$. This implies
$\pi(w\alpha)=\left(\begin{array}{cc}a&b\cr 0&a^{-1} \end{array}\right)$ and
$\pi(w')=\left(\begin{array}{cc}a'&b'\cr c'&d' \end{array}\right)$ with $d'\not=0$.
We get hence
$$\pi(w\alpha w')=\left(\begin{array}{cc}aa'+bc'&ab'+bd'\cr a^{-1}c'&a^{-1}d' 
\end{array}\right)$$
which shows $w\alpha w'\in {\cal C}$. The case $w\in {\cal C},\ w'\in {\cal A}$ 
follows using assertion (iv) of Theorem 2.2.
\medskip
\par {\bf Remark 3.3.} One can also use the more conceptual computation
$$\pi(w\alpha w')L(\infty)=\pi(w)\pi(\alpha)\pi(w')L(\infty)=
\pi(w)\pi(\alpha)L(0)=\pi(w)L(\infty)=L(0)$$
(for $w,w'\in{\cal A}$) in order to prove assertion (i) of Lemma 3.2.
\par Similarly, the case $w\in {\cal C},\ w'\in {\cal A}$ is dealt by
$$\pi(w\alpha w')L(\infty)=\pi(w)\pi(\alpha)\pi(w')L(\infty)=
\pi(w)\pi(\alpha)L(0)=\pi(w)L(\infty)\not=L(0)$$
and assertion (iv) of Theorem 2.2 completes the proof of assertion (ii) in Lemma 3.2.
\medskip
\par {\bf Proof of Theorem 2.6.} Assertion (i) follows easily from the previous 
lemma and the definition of ${\cal P}$.
\par Assertion (ii) follows from assertion (ii) of Theorem 2.2.
\medskip
\par {\bf Proof of Corollary 2.7.} Let $k$ be the finite field with
$q=p^d$ elements. An exercice using Theorem 2.6 shows that we have
$$\sharp({\cal A}^{l+1})=\sum_{k=0}^{[l/2]} {l-k\choose k} q^k(q-1)^{l-2k}
\ .$$
Corollary 2.3 establishes then the result if $x$ is a power of a prime number.
The proof follows now from the fact that both sides are polynomials in $x$.
\medskip
\par {\bf Proof of Theorem 2.9.} Follows from assertions (iii) and (ii) in Theorem 2.2.
\medskip
\par {\bf Proof of Theorem 2.10.} The elements 
$$\pi(\alpha_0\alpha_1\dots \alpha_{q-1}),\pi(\alpha_1\dots\alpha_{q-1}\alpha_0),\dots
,\pi(\alpha_{q-1}\alpha_0\dots\alpha_{q-2})\in{\rm SL}_2(k)$$
are all conjugate and have hence a common order $t$ which obviously works.
\medskip
\par The easy proof of Proposition 2.13 is left to the reader.
\medskip
\par {\bf Proof of Theorem 2.15.} This result follows readily from Theorem 2.9,
assertion (i) of Theorem 2.2 and Corollary 2.3.
\bigskip
\par I thank J.P. Allouche, P. de la Harpe and J. Helmstetter 
for useful comments.
\par I thank also an anonymous referee 
for the remark that the paper deals in fact with the projective group
${\rm PSL}_2(k)$ and for suggesting Remarks 3.1 and 3.3.
\bigskip
\par {\bf Bibliography}
\bigskip
\par {\bf [B1]} R. Bacher, {\it Curvature flow of maximal integral triangulations}, 
Ann. Inst. Fourier {\bf [49]}, 4 (1999), 1115-1128.
\par {\bf [B2]} R. Bacher, {\it An equivalence relation on $\{0,1\}^*$}, 
Europ. Journal of Combinatorics {\bf 21} (2000), 853-864.
\bigskip

Roland Bacher

INSTITUT FOURIER

Laboratoire de Math\'ematiques

UMR 5582 (UJF-CNRS)

BP 74

38402 St MARTIN D'H\`ERES Cedex (France)

e-mail: Roland.Bacher@ujf-grenoble.fr

\end{document}